\documentclass[12pt,english,a4paper]{article}
\usepackage[T1]{fontenc}
\usepackage[latin1]{inputenc}

\makeatletter

\newenvironment{lyxlist}[1]
{\begin{list}{}
{\settowidth{\labelwidth}{#1}
 \setlength{\leftmargin}{\labelwidth}
 \addtolength{\leftmargin}{\labelsep}
 }}
{\end{list}}


\usepackage{babel}
\makeatother

\begin{document}

\title{Mahlburg's work on Crank Functions returns to Ramanujan's work and
inspiration}

\author{Nagesh Juluru%
\thanks{Currently is a second year undergraduate student (B.Stat II Year)
at the Indian Statistical Institute. He has participated in this work
during August - September 2007 under the supervision of the second
author. %
} and Arni S.R. Srinivasa Rao %
\thanks{This work was initiated and continued when the author was at the Mathematical
Institute, University of Oxford, England, OX1 3LB. Email: arni@maths.ox.ac.uk. %
},%
\thanks{Corresponding author. Currently is Assistant Professor at the Indian
Statistical Institute. Address: 5th Floor, PSU, R.A. Fisher Building,
Indian Statistical Institute, 203 B.T. Road, Kolkata, INDIA 700108.
Email: arni@isical.ac.in%
} }

\maketitle
\begin{abstract}
Mahlburg (2005) brilliantly showed the importance of crank functions
in partition congruences that were originally guessed by Dyson (1944).
Ramanujan's partition functions are the centre of these works. Not
only for the theory on cranks, but for many other researchers' in
India Ramanujan's work inspired for their career in mathematics. This
is an undergraduate expository article. 
\end{abstract}

\section{Motivation }

Karl Mahlburg's Ph.D. dissertation work (2005) solidly established
the role of cranks (originally conceptualized by Freeman Dyson in
1944 and details are in section 3) in the theory of partition congruence
functions. His profound theory shows that indeed crank plays a central
role in explaining partitions congruences of Ramanujan-type and it
has once again oponed-up discussions on congruence modulo functions
by Srinivasa Ramanujan, a great mathematician of modern India. Mahlburg's
theory not only gave general ideas on crank functions but showed a
novel way to generate them through crank generating functions. He
proved Ken Ono's conjecture (2000), which has connections to the earlier
works by Andrews and Garvan (2000), Dyson (1944) and Ramanujan (1919).
Actual idea of crank is due to Dyson, while he was pursuing undergraduate
course at the University of Cambridge and wanted to explain Ramanujan's
congruence modulo 11. Ramanujan made substantial contributions in
the analytical theory of numbers, elliptic functions, continued functions
and infinite series. His conjectures in mathematics are still being
researched and new theories are being added. For a list of ongoing
research associated with Ramanujan's work and theories named after
him, see an article by Ken Ono in the \emph{Notices} (2006).

Proceedings of the National Academy of Sciences (PNAS) announced paper
of the year prize for the year 2005 to Karl Mahlburg, then a doctoral
candidate in mathematics at the University of Wisconsin, Madison.
Mahlburg's award winning paper in PNAS was entitled {}``Partition
congruence's and the Andrews-Garvan-Dyson crank.'' This article gave
a theoretical backing to the cranks and proved very interesting results
through crank generating functions. Ramanujan (1919) while studying
the number of distinct partitions of natural numbers $1$ to $200$
observed the following pattern

\begin{eqnarray*}
p(4),\: p(9),\: p(14),... & \cong & 0\:(mod\:5)\\
p(5),\: p(12),\: p(19),... & \cong & 0\:(mod\:7)\\
p(6),\: p(17),\: p(28),... & \cong & 0\:(mod\:11)\\
p(24),\: p(49),\: p(74),... & \cong & 0\:(mod\:25)\\
... & ... & ...\end{eqnarray*}

Using the idea of rank of the partition (see section 3), Dyson proposed
conjectures for the partitions $p(5n+4)$ and $p(7n+5).$ Then for
explaining the third conjecture using $p(11n+6)$, he guessed that
a crank of the partition would be helpful (see section 3). On the
other side, there were some attempts to improve the above partition
results for other primes after Ramanujan. Some how research in this
direction did not advance much for several decades (until the works
of Andrews, Garvan, Ahlgren, Ono) and researchers could obtain congruences
only for primes less than $33$. Mahlburg's outstanding paper improved
results that were given by Ono (2000) to prove the existence of infinite
families of partition congruences for the primes greater than $3.$
Before, we state Mahlburg's main results, we begin with some known
useful definitions in the theory numbers. Our focus of this work is
to understand Mahlburg's work on crank and give readers a short account
of the developments on crank to study Ramanujan-type functions and
also highlight how researchers in India inspired by Ramanujan's work
in general. In the next section, we introduce some elementary definitions
on partition functions, congruences, generating function and crank.
Section 3 begins from Dyson's guesses on crank's theory and a very
brief outline on this theory that was taken by Mahlburg to new heights.
In Section 4, we list selected recent outstanding results in the theory
of numbers in the world. In section 5, we mention selected important
contributions in the theory of numbers from India during post-Ramanujan
era (excluding the published work during last three decades) and author's
interactions with some of them.

\section{Useful definitions }

\subsection{Partition}

A partition of a positive number $n$ is any non-increasing sequence
of positive integers whose sum is $n.$ The partition function $p(n)$
is defined to count the number of distinct partitions of a given positive
integer $n.$ The rank of a partition is defined to be its largest
part minus the number of its parts. Examples: $(i)$ $p(8)=22$, a
partition of $8$ is $3+2+2+1$ and corresponding rank is $-1$. $(ii)$
$p(7)=14$, a partition of $7$ is $5+1+1$ and corresponding rank
is $2$ $(=5-3).$

\subsection{Congruence}

If $n$ is a positive integer, then two integers $a$ and $b$ are
called congruent modulo $n$, if $a-b$ is divisible by $n$ (they
give same remainder when divided by $n$). It is denoted by $a\cong b(mod\: n).$
For example $3462$ and $8649$ are congruent modulo $1729$ i.e.
$8649\cong3462(mod\:1729).$ Because difference between $3462$ and
$8649$ is $5187$ which is divisible by $1729.$ Both of these numbers
give same remainder when divided by $1729$ ($8649=1729\times5+4$
and $3462=1729\times2+4$).

Combining the idea of partitions and congruences, Ramanujan (1919)
for the first time introduced the idea of partition congruences appear
in section 1 and section 3.

\subsection{Crank}

Andrews and Garvan (2000) defined crank as follows: If $C=C_{1}+C_{2}+...+C_{s}+1+...+1$
has exactly $r$ $1s$, then let $o(C)$ be the number of parts of
$C$ that are strictly larger than $r$. The crank is given by

\begin{eqnarray*}
Crank(C) & = & \left\{ \begin{array}{cc}
C_{1} & \mbox{if}\: r=0\\
o(C)-r & \mbox{if}\: r\geq1\end{array}\right.\\
\\\end{eqnarray*}

For example, a partition of $7$ is $4+2+1$. Here, $o(7)=2$ and
$r=1.$ Crank of $7$ for this partition is$1.$ Similarly, a partition
of $9$ is $4+3+2$ and Crank of $9$ for this partition is $4.$

\subsection{Generating function}

It is a formal power series whose coefficients encode information
about a sequence $(a_{n}),$ where $n\in N.$ The generating function
of a sequence $(a_{n})$ is given by

\begin{eqnarray*}
F(a_{n};x) & = & \sum_{n\geq0}a_{n}x^{n}.\end{eqnarray*}

In a similar way, the generating function of a sequence $(a_{n_{1},n_{2},...,n_{k}})$,
where $n_{1},n_{2},...,n_{k}\in N$ is \begin{eqnarray*}
F(a_{n_{1},n_{2},...,n_{k}};x_{1},x_{2},...,x_{3}) & = & \sum_{n_{1}\geq0}\sum_{n_{2}\geq0}...\sum_{n_{k}\geq0}a_{n_{1},n_{2},...,n_{k}}x_{1}^{n_{1}}x_{2}^{n_{2}}...x_{k}^{n_{k}}.\end{eqnarray*}

\section{Evolving of the Mahlburg's work}

Mahlburg's outstanding work in 2005 on crank functions led to new
aspirations for further developments in the theory. He has shown that
his theory on crank functions satisfies wider range of congruences
than those predicted by earlier researchers like Ahlgren, Andrews,
Garvan, Dyson, and Ono. His approach of crank generating functions
together with Klein forms is novel. In this section, we begin with
a brief description of Ramanujan's congruences which actually led
Dyson to initiate work on cranks.

\subsection{Ramanujan's congruence functions and origins of Dyson's Crank }

Ramanujan's congruences for the partition function $p(n)$ in section
1 can be shortened into,

\begin{eqnarray*}
p(5n+4) & \cong & 0\:(mod\:5)\\
p(7n+5) & \cong & 0\:(mod\:7)\\
p(11n+6) & \cong & 0\:(mod\:11)\end{eqnarray*}

Above congruences can be written succinctly in the form $p(ln-\theta_{l})\cong0\:(mod\: l)$,
if we define $\theta_{l}=(l^{2}-1)/24$ (See Ahlgren and Ono, 2000).
Ramanujan proved first two congruences in his original paper (1919)
and later in 1921 (after his death in 1920), G.H. Hardy extracted
proofs of all the three congruences from an unpublished manuscript
of Ramanujan. Later, Watson (1938) published proofs for the congruences
after correcting them and also proved Ramanujan's congruence for the
modulus $11^{s}$ for some arbitrary $s.$ Lehner (1943, 1950) established
further properties of congruence modulo $11$. In 1944, Dyson published
a note {}``Some guesses in the Theory of Partitions'' where he has
given some combinatorial arguments to the above partition functions.
In this work he introduced notation $N(m,k)$ to denote the number
of partitions of $k$ with rank $m$ and $N(m,q,k)$ to denote the
number of partitions of $k$ whose rank is congruent to $m$ modulo
$q$. He obtained the rank of the partition by subtracting the number
of parts in a partition from the largest part. According to his example,
ranks of the partitions of $k$ are $k-1,k-3,k-4,...,2,1,0,-1,-2,...,4-k,3-k,1-k$
and

\begin{eqnarray*}
N(m,q,k) & = & \sum_{x=-\infty}^{\infty}N(m+xq,k).\end{eqnarray*}
 He conjectured that the partitions of $5k+4$ and $7k+5$ are divided
respectively into five and seven equal classes. Notationally these
conjectures are stated as

\begin{eqnarray*}
N(0,5,5k+4) & = & N(1,5,5k+4)\\
 &  & ,...,N(4,5,5k+4)\\
N(0,7,7k+5) & = & N(1,7,7k+5)\\
 &  & ,...,N(6,7,7k+5)\end{eqnarray*}
 By using the principle of conjugacy we can see that $N(m,k)=N(-m,k)$
and $N(m,q,k)=N(q-m,q,k)$. Dyson used this principle and reduced
the statements of above conjectures into smaller independent identities.
To simplify computing procedures Dyson translated these identities
into algebraic forms by the means of generating functions. Instead
of stating a similar conjecture for the partitions of $11k+6$, Dyson
predicted that conjecture with modulus $11$ is false and introduced
a hypothetical {}``crank'' of partition, $M(m,q,k)$, which explains
the number of partitions of $k$ whose crank is congruent to $m$
modulo $q$.

By using the relation $M(m,q,k)=M(q-m,q,k)$, Dyson guessed that

\begin{eqnarray*}
M(0,11,11k+6) & = & M(1,11,11k+6)\\
 &  & =M(2,11,11k+6)\\
 &  & \:\:\:=M(3,11,11k+6)\\
 &  & \:\:\:\:=M(4,11,11k+6)\end{eqnarray*}

\subsection{After Dyson's crank}

Atkin and Swinnerton-Dyer (1954) proved conjectures of Dyson and Atkin
extended Ramanujan's famous congruences to arbitrary powers of 5,
7, and 11 in his paper published in 1967. Working on Dyson's cranks,
Andrews and Garvan gave a formal definition to the crank (see section
2). Ken Ono has been researching on Ramanujan's contributions and
it is evident from his article in \emph{Notices} (2006) that he is
a fan of Ramanujan's work since his senior high school days. Ono (2000)
and Ahlgren and Ono (2001) advanced the idea on cranks and gave further
generalizations. Ono (2000) in his seminal work, established that
it is possible to express every prime $\geq5$ in the Ramanujan-type
congruences, followed by Ahlgren and Ono (2001), who laid complete
theoretical framework to describe such congruences.

\subsubsection*{Theorem (Ono):}

Let $l\geq5$ be prime and let $k$ be a positive integer. A positive
proportion of the primes $m$ have a property that

\begin{eqnarray*}
p\left(\frac{l^{k}m^{3}n+1}{24}\right) & \cong & 0\:(mod\: l)\end{eqnarray*}
 for every non negative integer $n$ coprime to $m.$

\subsubsection*{Theorem (Ahlgren and Ono)}

Define the integer $x_{l}\in\{\pm1\}$ by $x_{l}:=\left(\frac{-6}{l}\right)$
for each prime $l\geq5$ and let $S_{l}$ denote the set of $(l+1)/2$
integers $S_{l}:=\left\{ y\in(0,1,...,l-1):\left(\frac{y+\theta_{l}}{l}\right)=0\mbox{ or }-x_{l}\right\} ,$
where $\theta_{l}=(l^{2}-1)/24.$ Then Ahlgren and Ono proved the
monumental result for general partition functions: If $l\geq5$ is
prime, $k$ is a positive integer, and $y\in S_{l}$, then a positive
proportion of the primes $I\cong-1\:(mod\:24l)$ have the property
that

\begin{eqnarray*}
p\left(\frac{I^{3}n+1}{24}\right) & \cong & 0\;(mod\: l^{k})\end{eqnarray*}
 for all $n=1-24y\:(mod\:24l)$ with $(I,n)=1.$

Mahlburg (2005) in his award winning work improved Ono's results for
the primes $\geq5.$ By using crank generating functions (defined
by Andrews and Garvan (1988)) and elegantly linking them to Klein
forms he has demonstrated that the role of crank plays indeed a central
role in understanding the Ramanujan-type congruences. Statement of
his theorem is as follows: Let $l\geq5$ be primes and let $i$ and
$j$ be positive integers, then there are infinitely many arithmetic
progressions $Ak+B$ such that $M(m,l^{j},Ak+B)\cong0\:(mod\: l^{i})$
simultaneously for every $0\leq m\leq l^{j}-1.$

\section{Selected achievements in the long-standing problems in the theory
of numbers}

\begin{lyxlist}{00.00.0000}
\item [{{(i)}}] Fermat's last theorem was proved by Andrew Wiles (Princeton)
in 1995 . 
\item [{{(ii)}}] Manjul Bhargava (Princeton) generalized higher decomposition
laws in 2002 which were originally formulated by Gauss about 200 years
ago. 
\item [{{(iii)}}] Manindra Agrawal (Kanpur) along with his students Neeraj
Kayal and Nitin Saxena invented a computational algorithm in 2002
that can detect a given number is prime or not. 
\end{lyxlist}

\section{Selected contributions from India in the theory of numbers during
post-Ramanujan era}

After Ramanujan's work there were several researchers (inspired by
his legacy) from Indian universities and Institutes conducted serious
research in the theory of numbers. Some of the earliest researchers
in the theory of numbers after Ramanujan were S. Chowla, D.B. Lahiri,
M. Perisastri and J.M. Gandhi. See references for the list of their
selected publications. Interested readers can browse a larger publication
list of these authors from MathSciNet. We are not mentioning all the
contributions that came from India during the past three decades.

\begin{lyxlist}{00.00.0000}
\item [{{(i)}}] S. Chowla had wide range of interests from integers to
modular functions and contributed important results in the analytical
number theory. He was inspired by Ramanujan's mathematics and published
several papers in the theory of numbers. His personal biography is
found in the Internet. 
\item [{{(ii)}}] D. B. Lahiri extensively published on partition functions
and Ramanujan's congruences and other Ramanujan's functions and brought
some original thoughts to his work during 1960s. I%
\footnote{Arni S. R. Srinivasa Rao%
} have beautiful memories of him sharing experiences and listen to
his fond memories of early childhood days of mathematics in Burma
(now known as Myanmar) and how he got inspired by the works of Ramanujan.
He later moved to Kolkata to pursue his work in the theory of numbers.
He has well-known contributions in the sampling theory of statistics
and his method of sampling is frequently used in population surveys. 
\item [{{(iii)}}] M. Perisastri published novel results on Fermat's last
theorem and odd perfect numbers (using classical elementary arguments)
during 1950's and 1970's. Especially, his elementary results on Fermat's
last theorem were cited by many well known mathematicians working
on number theory in the world. M. Perisastri was a college level lecturer
with a depth of knowledge like a world class brilliant Professor in
mathematics. His inspiration by the work of Ramanujan was clear while
he was lecturing and through his beautiful naration of several incidences
between Hardy and Ramanujan. during Twice he was offered visiting
Professorship in a reputed US university, which he could not accept
due to family reasons. Many undergraduates from his home town Vizianagaram,
benefited by his teaching. He did not had much access to the publications
(except delayed arrival of \emph{Monthly, Crelle's Journal and other
American, British and Indian Journals}). He is known for composing
devotional songs and still he is active in singing songs in addition
to teaching undergraduate level mathematics. He has won several teaching
awards and wrote text books for college level calculus, number theory.
I%
\footnote{Arni S.R. Srinivasa Rao%
} have benefited by his inspirational mathematics teaching during my
high school and college level learning. 
\item [{{(iv)}}] J. M. Gandhi published several interesting results on
Fermat's last theorem and modular functions (during 1950s and 70s).
He was a very clear expositor in the theory of numbers and prolific
publisher of research papers. Most of his publications came while
he was working as a Professor at Jaipur and Bhilwara. 
\end{lyxlist}
We believe that Mahlburg's extraordinary paper has inspired many undergraduate
and graduate students. Some of them might look into Ramanujan's work
published almost a century ago and develop very good research questions.
We hope there will be further developments in crank functions by adopting
other kinds of generating functions and proving alternative results
to that of Mahlburg. We look for more papers that discover links between
Ramanujan's mathematics and present day science.

\end{document}